\documentclass[12pt,reqno]{amsart}
\usepackage{cite}
\usepackage{amsfonts}
\usepackage{amssymb, amsmath, latexsym}
\usepackage{longtable}
\usepackage{lipsum}
\usepackage{hyperref}
\usepackage{mleftright}
\usepackage{graphicx,mathtools}
\usepackage{mathtools}

\usepackage{array}

\numberwithin{equation}{section}
\theoremstyle{plain}
\newtheorem{theorem}{Theorem}[section]
\newtheorem{definition}{Definition}[section]

\title[Arithmetic properties of an analogue of $t$-core partitions] {Arithmetic properties of an analogue of $t$-core partitions}
\author[P. Talukdar]{Pranjal Talukdar}
\address{Department of Mathematical Sciences, Tezpur University,  Assam 784028, India}
\email{pranjaltalukdar113@gmail.com}

\allowdisplaybreaks
\begin{document}
\begin{abstract}
An integer partition of a positive integer $n$ is called to be $t$-core if none of its hook lengths are divisible by $t$.  Recently, Gireesh, Ray and Shivashankar [`A new analogue of $t$-core partitions', \textit{Acta Arith.} \textbf{199} (2021), 33-53] introduced an analogue $\overline{a}_t(n)$ of the $t$-core partition function $c_t(n)$. They obtained certain multiplicative formulas and arithmetic identities for $\overline{a}_t(n)$ where $t \in \{3,4,5,8\}$ and studied the arithmetic density of $\overline{a}_t(n)$ modulo $p_i^{j}$ where $t=p_1^{a_1}\cdots p_m^{a_m}$ and $p_i\geq 5$ are primes. Very recently, Bandyopadhyay and Baruah [`Arithmetic identities for some analogs of the
5-core partition function', \textit{J. Integer Seq.} \textbf{27} (2024), \# 24.4.5] proved new arithmetic identities satisfied by $\overline{a}_5(n)$. In this article, we study the arithmetic densities of $\overline{a}_t(n)$ modulo arbitrary powers of 2 and 3 for $t=3^\alpha m$ where $\gcd(m,6)$=1. Also, employing a result of Ono and Taguchi on the nilpotency of Hecke operators, we prove an infinite family of congruences for $\overline{a}_3(n)$ modulo arbitrary powers of 2.

\end{abstract}
\maketitle
\noindent{\footnotesize Key words: $t$-core partition, analogue of $t$-core partition, theta functions, modular forms, arithmetic density}

\vskip 3mm
\noindent {\footnotesize 2010 Mathematical Reviews Classification Numbers: 11P83, 05A17, 11F11}

\section{Introduction and statement of results}\label{Introduction}

A partition $\pi=\{\pi_1,\pi_2,\cdots,\pi_k\}$ of a positive integer $n$ is a non-increasing sequence of natural numbers such that $\displaystyle\sum_{i=1}^k \pi_i=n$. The number of partitions of $n$ is denoted by $p(n)$. The Ferrers–Young diagram of $\pi$ is an array of nodes with $\pi_i$ nodes in the $i$th row. The $(i,j)$ hook is the set of nodes directly to the right of $(i,j)$ node, together with the set of nodes directly below it, as well as the $(i,j)$ node itself. The hook number, $H(i,j)$, is the total number of nodes on the $(i,j)$ hook. For a positive integer $t\geq 2$, a partition of $n$ is said to be $t$-core if none of the hook numbers are divisible by $t$. We illustrate the Ferrers-Young diagram of the partition $4+3+1$ of 8 with hook numbers as follows:
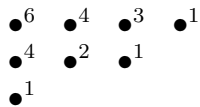
\begin{center}
\begin{tabular}{ cccc }  
 $\bullet^6$ &  $\bullet ^4$&  $\bullet^3$&  $\bullet^1$\\ 
 $ \bullet^4$& $ \bullet^2$ &  $\bullet^1$&\\ 
 $\bullet^1$ &  & & \\ 
\end{tabular}
\end{center}
It is clear that for $t\geq7$, the partition $4+3+1$ of 8 is a $t$-core partition.

Suppose that $c_t(n)$ counts the $t$-core partitions of $n$, then the generating function of $c_t(n)$ is given by (see \cite[Eq 2.1]{gks})
\begin{align}
\displaystyle\sum_{n=0}^{\infty}c_t(n)q^n=\dfrac{(q^t;q^t)_\infty^t}{(q;q)_\infty}=\dfrac{f_t^t}{f_1}\label{c t},
\end{align}
where here and throughout the paper, for $|q| <1$, we define $(a;q)_\infty:=\displaystyle\prod_{k=0}^\infty(1-aq^k)$ and for convenience, we set $f_k:=(q^k;q^k)_\infty$ for integers $k\geq 1$.

In an existence result, Granville and Ono \cite{go} proved that if $t\geq 4$, then $c_t(n)>0$ for every nonnegative integer $n$. A brief survey of $t$-core partitions can be found in \cite{ck}.

Again, for an integral power series $F(q):=\displaystyle\sum_{n=0}^{\infty}a(n)q^n$ and $0\leq r < M$, we define the arithmetic density $\delta_r(F,M;X)$ as
\begin{align*}
\delta_r(F,M;X):=\dfrac{\#\left\{0\leq n\leq X: a(n)\equiv r\pmod M\right\}}{X}.
\end{align*}
An integral power series $F$ is called \textit{lacunary modulo} $M$ if
\begin{align*}
\lim_{X\rightarrow\infty}\delta_0(F,M;X)=1,
\end{align*}
that is, almost all of the coefficients of $F$ are divisible by $M$.

Arithmetic densities of $c_t(n)$ modulo arbitrary powers of 2, 3 and primes greater than or equal to 5 are recently studied by Jindal and Meher \cite{jm}. 

Recall that for $\vert ab\vert <1$, Ramanujan's general theta function $f(a,b)$ is given by
\begin{align*}
f(a,b)=\displaystyle\sum_{n=-\infty}^\infty a^{n(n+1)/2}b^{n(n-1)/2}.
\end{align*}
In Ramanujan's notation, the Jacobi triple product identity \cite[p. 35, Entry 19]{bcb3} takes the shape
\begin{align*}
f(a,b)=(-a;ab)_\infty(-b;ab)_\infty(ab;ab)_\infty.
\end{align*}

Consider the following two special cases of $f(a,b)$:
\begin{align}
\varphi(-q)&:=f(-q,-q)=\displaystyle\sum_{n=-\infty}^\infty (-1)^nq^{n^2}=\dfrac{f_1^2}{f_2}\label{p-q},\\
f(-q)&:=f(-q,-q^2)=\displaystyle\sum_{n=-\infty}^\infty (-1)^nq^{n(3n-1)/2}=f_1\label{f-q}.
\end{align}

In the notation of Ramanujan's theta functions, the generating function of $c_t(n)$ may be rewritten as
\begin{align}
\displaystyle\sum_{n=0}^{\infty}c_t(n)q^n=\dfrac{f^t(-q^t)}{f(-q)}.\label{ct1}
\end{align}

Recently, Gireesh, Ray and Shivashankar \cite{grs} considered a new function $\overline{a}_t(n)$ by substituting $\varphi(-q)$ in place of $f(-q)$ in the generating function of $c_t(n)$ (in \eqref{ct1}), namely
\begin{align}
\displaystyle\sum_{n=0}^{\infty}\overline{a}_t(n)q^n=\dfrac{\varphi^t(-q^t)}{\varphi(-q)}=\dfrac{f_2f_t^{2t}}{f_1^2f_{2t}^t}\label{a t}.
\end{align}
They proved several multiplicative formulas and arithmetic identities for $\overline{a}_t(n)$ for $t=$ 2, 3, 4 and 8 using Ramanujan's theta functions and $q$-series techniques. Using the theory of modular forms, they studied the divisibility of $\overline{a}_t(n)$ modulo arbitrary powers of primes greater than 5. More precisely, they proved the following theorem.
\begin{theorem}
Let $t=p_1^{a_1}\cdots p_m^{a_m}$ where $p_i$’s are prime numbers greater than or equal to 5. Then for every positive integer $j$, we have
\begin{align*}
\lim_{X\rightarrow\infty}\dfrac{\#\left\{0\leq n\leq X: \overline{a}_t(n)\equiv 0\pmod{p_i^j}\right\}}{X}=1.
\end{align*}
\end{theorem}
They also deduced a Ramanujan type congruence for $\overline{a}_5(n)$ modulo 5 by using an algorithm developed by Radu and Sellers \cite{rs}. Bandyopadhyay and Baruah \cite{bb} proved some new identities connecting $\overline{a}_5(n)$ and $c_5(n)$. They also found a reccurence relation for $\overline{a}_5(n)$.

Recently, Cotron \textit{et. al.} \cite[Theorem 1.1]{cot} proved a strong result regarding lacunarity of eta-quotients modulo arbitrary powers of primes under certain conditions. We observe that the eta-quotients associated with $\overline{a}_t(n)$ do not satisfy these conditions, which makes the problem of studying lacunarity of $\overline{a}_t(n)$ more interesting. In this article, we study the arithmetic densities of $\overline{a}_{t}(n)$ modulo arbitrary powers of 2 and 3 where $t=3^\alpha m$. To be specific, we prove the following theorems.

\begin{theorem}\label{a 2k}
Let $k\geq 1$, $\alpha\geq 0$ and $m \geq 1$ be integers with \textup{gcd}$(m, 6)=1$. Then the
set
\begin{align*}
\left\{n \in \mathbb{N} : \overline{a}_{3^\alpha m}(n) \equiv 0 \pmod {2^k}\right\}
\end{align*}
has arithmetic density 1.
\end{theorem}

\begin{theorem}\label{a 3k}
Let $k\geq 1$, $\alpha\geq 0$ and $m \geq 1$ be integers with \textup{gcd}$(m, 6)=1$. Then the
set
\begin{align*}
\left\{n \in \mathbb{N} : \overline{a}_{3^\alpha m}(n) \equiv 0 \pmod {3^k}\right\}
\end{align*}
has arithmetic density 1.
\end{theorem}

The fact that the action of Hecke algebras on spaces of modular forms of level 1 modulo 2 is locally nilpotent was first observed by Serre and proved by Tate (see \cite{ser1}, \cite{ser2}, \cite{ta}). Later, this result was generalized to higher levels by Ono and Taguchi \cite{ot}. In this article, we observe that the
eta-quotient associated to $\overline{a}_3(n)$ is a modular form whose level is in
the list of Ono and Taguchi. Thus, we use a result of Ono and Taguchi to prove the following congruences for $\overline{a}_3(n)$.

\begin{theorem}\label{cong1}
Let $n$ be a nonnegative integer. Then there exists an integer $c\geq 0$
such that for every $d\geq 1$ and distinct primes $p_1,\dots,p_{c+d}$ coprime to 6, we have
\begin{align*}
\overline{a}_3\left(\dfrac{p_1\cdots p_{c+d}\cdot n}{24}\right)\equiv 0\pmod{2^d}
\end{align*}
whenever $n$ is coprime to $p_1,\dots,p_{c+d}$.
\end{theorem}

The paper is organized as follows. In Section \ref{prelim}, we state some preliminaries of the theory of modular forms.   Then we prove Theorems \ref{a 2k}-\ref{cong1} using the properties of modular forms in Sections \ref{pr a 2k}-\ref{pr cong1} respectively. And finally, we mention some directions for future study in the concluding section. 
\section{Preliminaries}\label{prelim}

In this section, we recall some basic facts and definitions on modular forms. For more details, one can see \cite{kob} and \cite{ono}.

Firstly, we define the matrix groups
\begin{align*}
\textup{SL}_2(\mathbb{Z})&:=\left\{\begin{bmatrix}
a & b \\
c & d 
\end{bmatrix} : a, b, c, d \in \mathbb{Z}, ad-bc=1\right\},\\
\Gamma_0(N)&:=\left\{\begin{bmatrix}
a & b \\
c & d 
\end{bmatrix} \in \textup{SL}_2(\mathbb{Z}) : c\equiv 0 \pmod N\right\},\\
\Gamma_1(N)&:=\left\{\begin{bmatrix}
a & b \\
c & d 
\end{bmatrix} \in 
\Gamma_0(N) : a\equiv d\equiv 1 \pmod N\right\},
\end{align*}
and
\begin{align*}
\Gamma(N):=\left\{\begin{bmatrix}
a & b \\
c & d 
\end{bmatrix} \in 
\textup{SL}_2(\mathbb{Z}) : a\equiv d\equiv 1 \pmod N, \textup{ and } b\equiv c\equiv 0 \pmod N \right\},
\end{align*}
where $N$ is a positive integer. A subgroup $\Gamma$ of $\textup{SL}_2(\mathbb{Z})$ is called a congruence subgroup if $\Gamma(N)\subseteq \Gamma$ for some $N$ and the smallest $N$ with this property is called the level of $\Gamma$. For instance, $\Gamma_0(N)$ and $\Gamma_1(N)$ are congruence subgroups of level $N$.

Let $\mathbb{H}$ denote the upper half of the complex plane. The group 
\begin{align*}
\textup{GL}_2^{+}(\mathbb{R}):=\left\{\begin{bmatrix}
a & b \\
c & d 
\end{bmatrix} \in 
\textup{SL}_2(\mathbb{Z}) : a, b, c, d \in \mathbb{R} \textup{ and } ad-bc>0 \right\}
\end{align*}
acts on $\mathbb{H}$ by $\begin{bmatrix}
a & b \\
c & d 
\end{bmatrix} z = \dfrac{az+b}{cz+d}$. We identify $\infty$ with $\dfrac{1}{0}$ and define $\begin{bmatrix}
a & b \\
c & d 
\end{bmatrix} \dfrac{r}{s} = \dfrac{ar+bs}{cr+ds}$, where $\dfrac{r}{s} \in \mathbb{Q} \cup \{ \infty \}$. This gives an action of $\textup{GL}_2^{+}(\mathbb{R})$ on the extended upper half plane $\mathbb{H}^*=\mathbb{H}\cup \mathbb{Q}\cup \{\infty\}$. Suppose that $\Gamma$ is a congruence subgroup of $\textup{SL}_2(\mathbb{Z})$. A cusp of $\Gamma$ is an equivalence class in $\mathbb{P}^1=\mathbb{Q}\cup\{\infty\}$ under the action of $\Gamma$. 

The group $\textup{GL}_2^{+}(\mathbb{R})$ also acts on functions $f:\mathbb{H}\rightarrow \mathbb{C}$. In particular, suppose that $\gamma=\begin{bmatrix}
a & b \\
c & d 
\end{bmatrix}$ $\in \textup{GL}_2^{+}(\mathbb{R})$. If $f(z)$ is a meromorphic function on $\mathbb{H}$ and $\ell$ is an integer, then define the slash operator $|_{\ell}$ by
\begin{align*}
(f|_{\ell}\gamma)(z):=(\textup{det}(\gamma))^{\ell/2}(cz+d)^{-\ell}f(\gamma z).
\end{align*}

\begin{definition}
Let $\Gamma$ be a congruence subgroup of level $N$. A holomorphic function $f:\mathbb{H}\rightarrow \mathbb{C}$ is called a modular form with integer weight $\ell$ on $\Gamma$ if the following hold:
\begin{enumerate}
\item[(1)] We have
\begin{align*}
f\left(\dfrac{az+b}{cz+d}\right)=(cz+d)^{\ell}f(z)
\end{align*}
for all $z \in \mathbb{H}$ and and all $\begin{bmatrix}
a & b \\
c & d 
\end{bmatrix}$ $\in \Gamma$.
\item[(2)] If $\gamma \in$ $\textup{SL}_2(\mathbb{Z})$, then $(f|_{\ell}\gamma)(z)$ has a Fourier expansion of the form 
\begin{align*}
(f|_{\ell}\gamma)(z)=\displaystyle\sum_{n\geq 0}a_{\gamma}(n)q_{N}^n,
\end{align*}
where $q:=e^{2\pi iz/N}$.
\end{enumerate}
\end{definition}
For a positive integer $\ell$, the complex vector space of modular forms of weight $\ell$ with respect to a congruence subgroup $\Gamma$ is denoted by $M_{\ell}(\Gamma)$.
\begin{definition}\cite[Definition 1.15]{ono}
If $\chi$ is a Dirichlet character modulo $N$, then we say that a modular form $f \in M_{\ell}(\Gamma_{1}(N))$ has Nebentypus character $\chi$ if 
\begin{align*}
f\left(\dfrac{az+b}{cz+d}\right)=\chi(d)(cz+d)^{\ell}f(z)
\end{align*}
for all $z\in \mathbb{H}$ and all $\begin{bmatrix}
a & b \\
c & d 
\end{bmatrix}$ $\in \Gamma_{0}(N)$. The space of such modular forms is denoted by $M_{\ell}(\Gamma_{0}(N),\chi)$.
\end{definition}

The relevant modular forms for the results of this paper arise from eta-quotients. Recall that the Dedekind eta-function $\eta(z)$ is defined by
\begin{align*}
\eta(z):=q^{1/24}(q;q)_{\infty}=q^{1/24}\displaystyle\prod_{n=1}^{\infty}(1-q^n),
\end{align*}
where $q:=e^{2\pi iz}$ and $z\in \mathbb{H}$. A function $f(z)$ is called an eta-quotient if it is of the form 
\begin{align*}
f(z)=\displaystyle\prod_{\delta|N}\eta(\delta z)^{r_{\delta}},
\end{align*}
where $N$ is a positive integer and $r_{\delta}$ is an integer. Now, we recall two important theorems from \cite[p. 18]{ono} which
will be used later.

\begin{theorem}\cite[Theorem 1.64]{ono}\label{level}
If $f(z)=\displaystyle\prod_{\delta|N}\eta(\delta z)^{r_{\delta}}$ is an eta-quotient such that $\ell=\dfrac{1}{2}\displaystyle\sum_{\delta|N}r_{\delta}$ $\in \mathbb{Z}$,
\begin{align*}
\displaystyle\sum_{\delta|N}\delta r_{\delta}\equiv 0 \pmod {24} \quad\quad \text{and}\quad\quad
\displaystyle\sum_{\delta|N}\dfrac{N}{\delta} r_{\delta}\equiv 0 \pmod {24},
\end{align*}
then $f(z)$ satisfies 
\begin{align*}
f\left(\dfrac{az+b}{cz+d}\right)=\chi(d)(cz+d)^{\ell}f(z)
\end{align*}
for every $\begin{bmatrix}
a & b \\
c & d 
\end{bmatrix}$ $\in \Gamma_{0}(N)$. Here the character $\chi$ is defined by $\chi(d)$ $:=$ $\left(\dfrac{(-1)^{\ell}s}{d}\right),$ where $s:=$ $\displaystyle\prod_{\delta|N}\delta^{r_{\delta}}$.
\end{theorem}

Consider $f$ to be an eta-quotient which satisfies the conditions of Theorem \ref{level} and that the associated weight $\ell$ is a positive integer. If $f(z)$ is holomorphic at all the cusps of $\Gamma_{0}(N)$, then $f(z) \in M_{\ell}\left(\Gamma_{0}(N), \chi\right)$. The necessary criterion for determining orders of an
eta-quotient at cusps is given by the following theorem.
\begin{theorem}\cite[Theorem 1.64]{ono}\label{order}
Let $c$, $d$ and $N$ be positive integers with $d|N$ and gcd$(c,d)$=1. If $f$ is an eta-quotient satisfying the conditions of Theorem \ref{level} for $N$, then the order of vanishing of $f(z)$ at the cusp $(c/d)$ is 
\begin{align*}
\dfrac{N}{24}\displaystyle\sum_{\delta|N}\dfrac{\text{gcd}(d,\delta)^2r_{\delta}}{\text{gcd}(d,N/d)d\delta}.
\end{align*}
\end{theorem}

We now recall a deep theorem of Serre \cite[Page 43]{ono} which will be used in proving Theorems \ref{a 2k} and \ref{a 3k}.
\begin{theorem} \cite[p. 43]{ono}\label{Ser}
Let $g(z) \in M_k(\Gamma_0(N),\chi)$ has Fourier expansion
\begin{align*}
g(z)=\displaystyle\sum_{n=0}^\infty b(n)q^n \in \mathbb{Z}[[q]].
\end{align*}
Then for a positive integer $r$, there is a constant $\alpha > 0$ such that
\begin{align*}
\#\{0<n\leq X: b(n)\not\equiv 0 \pmod r\}=\mathcal{O}\left(\dfrac{X}{(log X)^\alpha}\right).
\end{align*}
Equivalently
\begin{align*}
\displaystyle{\lim_{X \to \infty}} \dfrac{\#\{0<n\leq X: b(n)\not\equiv 0 \pmod r\}}{X}=0.
\end{align*}
\end{theorem}

Finally, we recall the definition of Hecke operators.  Let $m$ be a positive integer and $f(z)=\displaystyle\sum_{n=0}^{\infty}a(n)q^n$ $\in M_{\ell}(\Gamma_{0}(N),\chi)$. Then the action of Hecke operator $T_{m}$ on $f(z)$ is defined by
\begin{align*}
f(z)|T_m:=\displaystyle\sum_{n=0}^{\infty}\left(\displaystyle\sum_{d|\text{gcd}(n,m)}\chi(d)d^{\ell-1}a\left(\dfrac{nm}{d^2}\right)\right)q^n.
\end{align*}
In particular, if $m=p$ is prime, then we have
\begin{align}
\label{operator}f(z)|T_p:=\displaystyle\sum_{n=0}^{\infty}\left(a(pn)+\chi(p)p^{\ell-1}a\left(\dfrac{n}{p}\right)\right)q^n.
\end{align}
We note that $a(n)=0$ unless $n$ is a nonnegative integer.

\section{Proof of Theorem \ref{a 2k}}\label{pr a 2k}
Putting $t=3^{\alpha}m$ in \eqref{a t}, we have
\begin{align}
\displaystyle\sum_{n=0}^{\infty}\overline{a}_{3^{\alpha}m}(n)q^n=\dfrac{f_2f_{3^{\alpha}m}^{2\cdot3^{\alpha}m}}{f_1^2f_{2\cdot3^{\alpha}m}^{3^{\alpha}m}}\label{a 3am}.
\end{align}
We define
\begin{align*}
A_{\alpha,m}(z):=\dfrac{\eta^2\left(2^33^{\alpha+1}mz\right)}{\eta\left(2^43^{\alpha+1}mz\right)}.
\end{align*}
For any prime $p$ and positive integer $j$, we have
\begin{align*}
(q;q)_{\infty}^{p^j}\equiv (q^p;q^p)_{\infty}^{p^{j-1}}\pmod{p^j}.
\end{align*}
Using the above relation, for any integer $k\geq 1$, we get
\begin{align}
A^{2^k}_{\alpha,m}(z)=\dfrac{\eta^{2^{k+1}}\left(2^33^{\alpha+1}mz\right)}{\eta^{2^k}\left(2^43^{\alpha+1}mz\right)}\equiv 1\pmod{2^{k+1}}.\label{A 2k}
\end{align}
Next we define 
\begin{align*}
B_{\alpha,m,k}(z)&:=\dfrac{\eta(48z)\eta^{2\cdot3^{\alpha}m}\left(2^33^{\alpha+1}mz\right)}{\eta^2(24z)\eta^{3^{\alpha}m}\left(2^43^{\alpha+1}mz\right)}A^{2^k}_{\alpha,m}(z)\\
&=\dfrac{\eta(48z)\eta^{2\cdot3^{\alpha}m+2^{k+1}}\left(2^33^{\alpha+1}mz\right)}{\eta^2(24z)\eta^{3^{\alpha}m+2^k}\left(2^43^{\alpha+1}mz\right)}.
\end{align*}
In view of \eqref{a 3am} and \eqref{A 2k}, we have
\begin{align}
B_{\alpha,m,k}(z)&\equiv\dfrac{\eta(48z)\eta^{2\cdot3^{\alpha}m}\left(2^33^{\alpha+1}mz\right)}{\eta^2(24z)\eta^{3^{\alpha}m}\left(2^43^{\alpha+1}mz\right)}\notag\\
&\equiv\dfrac{f_{48}f_{2^33^{\alpha+1}m}^{2\cdot3^{\alpha}m}}{f_{24}^2f_{2^4\cdot3^{\alpha+1}m}^{3^{\alpha}m}}
\equiv\displaystyle\sum_{n=0}^{\infty}\overline{a}_{3^{\alpha}m}(n)q^{24n}\pmod{2^{k+1}}\label{B amk}.
\end{align}

Next, we will show that $B_{\alpha,m,k}(z)$ is a modular form. Applying Theorem \ref{level}, we find that the level of $B_{\alpha,m,k}(z)$ is $N=2^43^{\alpha+1}mM$, where $M$ is the smallest positive integer such that
\begin{align*}
2^43^{\alpha+1}mM\left(\dfrac{-2}{24}+\dfrac{1}{48}+\dfrac{2\cdot3^{\alpha}m+2^{k+1}}{2^33^{\alpha+1}m}+\dfrac{-3^{\alpha}m-2^k}{2^43^{\alpha+1}m}\right)\equiv 0\pmod{24},
\end{align*}
which implies 
\begin{align*}
3\cdot2^kM\equiv 0\pmod{24}.
\end{align*}
Therefore $M=4$ and the level of $B_{\alpha,m,k}(z)$ is $N=2^6 3^{\alpha+1}m$.

The representatives for the cusps of $\Gamma_0\left(2^63^{\alpha+1}m\right)$ are given by fractions $c/d$ where $d|2^63^{\alpha+1}m$ and $\gcd(c, 2^63^{\alpha+1}m) = 1$ (see \cite[Proposition 2.1]{cot}). By Theorem \ref{order}, $B_{\alpha,m,k}(z)$ is holomorphic at a cusp $c/d$ if and only if
\begin{align*}
-2 \dfrac{\gcd (d,24)^2}{24}+ \dfrac{\gcd (d,48)^2}{48}+\left(3^\alpha m+2^k\right)\left(2\dfrac{ \gcd \left(d,2^3 3^{\alpha+1} m\right)^2}{2^3 3^{\alpha+1} m}-\dfrac{\gcd \left(d,2^4 3^{\alpha+1} m\right)^2}{2^4 3^{\alpha+1} m}\right)\geq 0.
\end{align*}
Equivalently, $B_{\alpha,m,k}(z)$ is holomorphic at a cusp $c/d$ if and only if
\begin{align*}
L:=3^{\alpha}m(-4G_1+G_2+4G_3-1)+2^k(4G_3-1)\geq 0,
\end{align*}
where $G_1=\dfrac{\gcd (d,24)^2}{\gcd \left(d,2^4 3^{\alpha+1} m\right)^2}$, $G_2=\dfrac{\gcd (d,48)^2}{\gcd \left(d,2^43^{\alpha+1} m\right)^2}$  and 
$G_3=\dfrac{\gcd (d,2^3 3^{\alpha+1} m)^2}{\gcd \left(d,2^4 3^{\alpha+1} m\right)^2}$.

Let $d$ be a divisor of $2^6 3^{\alpha+1} m$. We
can write $d = 2^{r_1}3^{r_2}t$ where $0 \leq r_1 \leq 6$, $0 \leq r_2 \leq \alpha + 1$ and $t|m$. We now consider the
following two cases depending on $r_1$.

Case 1: Let $0\leq r_1\leq 3$, $0\leq r_2\leq \alpha+1$. Then $G_1=G_2$, $\dfrac{1}{3^{2\alpha }t^2}\leq G_1\leq 1$ and $G_3=1$. Therefore $L= 3^{\alpha+1}m(1-G_1)+3\cdot2^k\geq 3\cdot2^k.$

Case 2: Let $4\leq r_1\leq 6$, $0\leq r_2\leq \alpha+1$. Then $G_2=4G_1$, $\dfrac{1}{4\cdot3^{2\alpha }t^2}\leq G_1\leq \dfrac{1}{4}$ and $G_3=\dfrac{1}{4}$ which implies $L= 0.$

Hence, $B_{\alpha,m,k}(z)$ is holomorphic at every cusp $c/d$. The weight of $B_{\alpha,m,k}(z)$ is $\ell=\dfrac{1}{2} \left(3^\alpha m+2^k-1\right)$ which is a positive integer and the associated character is given by
 \begin{align*}
\chi_1(\bullet)=\left(\dfrac{(-1)^\ell 3^{(\alpha+1) \left(3^\alpha m+2^k\right)-1} m^{3^\alpha m+2^k}}{\bullet}\right).
\end{align*}

 Thus, $B_{\alpha,m,k}(z)\in M_\ell\left(\Gamma_0(N),\chi\right)$ where $\ell$, $N$ and $\chi$ are as above. Therefore by Theorem \ref{Ser}, the Fourier coefficients of $B_{\alpha,m,k}(z)$ are almost divisible by $r=2^k$. Due to \eqref{B amk}, this holds for $\overline{a}_{3^{\alpha}m}(n)$ also. This completes the proof of Theorem \ref{a 2k}.
\section{Proof of Theorem \ref{a 3k}}\label{pr a3k}
We proceed along the same lines as in the proof of Theorem \ref{a 2k}. Here we define
\begin{align*}
C_{\alpha,m}(z):=\dfrac{\eta^3\left(2^43^{\alpha+1}mz\right)}{\eta\left(2^43^{\alpha+2}mz\right)}.
\end{align*}
Using the binomial theorem, for any integer $k\geq 1$, we have
\begin{align}
C^{3^k}_{\alpha,m}(z)=\dfrac{\eta^{3^{k+1}}\left(2^43^{\alpha+1}mz\right)}{\eta^{3^k}\left(2^43^{\alpha+2}mz\right)}\equiv 1\pmod{3^{k+1}}\label{C 3k}.
\end{align}
Next we define 
\begin{align*}
D_{\alpha,m,k}(z)&:=\dfrac{\eta(48z)\eta^{2\cdot3^{\alpha}m}\left(2^33^{\alpha+1}mz\right)}{\eta^2(24z)\eta^{3^{\alpha}m}\left(2^43^{\alpha+1}mz\right)}C^{3^k}_{\alpha,m}(z)\\
&=\dfrac{\eta(48z)\eta^{2\cdot3^{\alpha}m}\left(2^33^{\alpha+1}mz\right)\eta^{3^{k+1}-3^{\alpha}m}\left(2^43^{\alpha+1}mz\right)}{\eta^2(24z)\eta^{3^{k}}\left(2^43^{\alpha+2}mz\right)}.
\end{align*}
From \eqref{a 3am} and \eqref{C 3k}, we have
\begin{align}
D_{\alpha,m,k}(z)&\equiv\dfrac{\eta(48z)\eta^{2\cdot3^{\alpha}m}\left(2^33^{\alpha+1}mz\right)}{\eta^2(24z)\eta^{3^{\alpha}m}\left(2^43^{\alpha+1}mz\right)}\notag\\
&\equiv\dfrac{f_{48}f_{2^33^{\alpha+1}m}^{2\cdot3^{\alpha}m}}{f_{24}^2f_{2^4\cdot3^{\alpha+1}m}^{3^{\alpha}m}}
\equiv\displaystyle\sum_{n=0}^{\infty}\overline{a}_{3^{\alpha}m}(n)q^{24n}\pmod{3^{k+1}}\label{D amk}.
\end{align}

We now prove that $D_{\alpha,m,k}(z)$ is a modular form. Applying Theorem \ref{level}, we find that the level of $D_{\alpha,m,k}(z)$ is $N=2^43^{\alpha+2}mM$, where $M$ is the smallest positive integer such that
\begin{align*}
2^43^{\alpha+2}mM\left(\dfrac{-2}{24}+\dfrac{1}{48}+\dfrac{2\cdot3^{\alpha}m}{2^33^{\alpha+1}m}+\dfrac{3^{k+1}-3^{\alpha}m}{2^43^{\alpha+1}m}+\dfrac{-3^k}{2^43^{\alpha+2}m}\right)\equiv 0\pmod{24},
\end{align*}
which gives 
\begin{align*}
8\cdot3^kM\equiv 0\pmod{24}.
\end{align*}
Therefore $M=1$ and the level of $D_{\alpha,m,k}(z)$ is $N=2^4 3^{\alpha+2}m$.

The representatives for the cusps of $\Gamma_0\left(2^43^{\alpha+2}m\right)$ are given by fractions $c/d$ where $d|2^43^{\alpha+2}m$ and $\gcd(c, 2^43^{\alpha+2}m) = 1$. By using Theorem \ref{order}, $D_{\alpha,m,k}(z)$ is holomorphic at a cusp $c/d$ if and only if
\begin{align*}
&-2 \dfrac{\gcd (d,24)^2}{24}+ \dfrac{\gcd (d,48)^2}{48}+2\cdot 3^{\alpha} m \dfrac{\gcd \left(d,2^33^{\alpha+1} m\right)^2}{2^33^{\alpha+1} m}+\left(3^{k+1}-3^{\alpha} m\right) \dfrac{\gcd \left(d,2^4 3^{\alpha+1} m\right)^2}{2^43^{\alpha+1} m}\\
&\quad\quad\quad\quad\quad\quad\quad\quad\quad\quad\quad\quad-3^k~ \dfrac{\gcd \left(d,2^4 3^{\alpha+2} m\right)^2}{2^43^{\alpha+2} m}
\geq 0.
\end{align*}
Equivalently, $D_{\alpha,m,k}(z)$ is holomorphic at a cusp $c/d$ if and only if
\begin{align*}
L:=3^{\alpha+1}m\left(-4G_1+G_2+4G_3-G_4\right)+3^{k}(9G_4-1)\geq 0,
\end{align*}
where $G_1=\dfrac{\gcd (d,24)^2}{\gcd \left(d,2^4 3^{\alpha+2} m\right)^2}$, $G_2=\dfrac{\gcd (d,48)^2}{\gcd \left(d,2^4 3^{\alpha+2} m\right)^2}$, $G_3=\dfrac{\gcd (d,2^3 3^{\alpha+1} m)^2}{\gcd \left(d,2^4 3^{\alpha+2} m\right)^2}$

 and $G_4=\dfrac{\gcd (d,2^4 3^{\alpha+1} m)^2}{\gcd \left(d,2^4 3^{\alpha+2} m\right)^2}$.

Let $d$ be a divisor of $2^43^{\alpha+2} m$. We write $d = 2^{r_1}3^{r_2}t$ where $0 \leq r_1 \leq 4$, $0 \leq r_2 \leq \alpha + 2$ and $t|m$. We now consider the
following four cases depending on the values of $r_1$ and $r_2$.

Case 1: Let $0\leq r_1\leq 3$, $0\leq r_2\leq \alpha+1$. Then $G_1=G_2$, $\dfrac{1}{3^{2\alpha}t^2}\leq G_1\leq 1$ and $G_3=G_4=1$. Hence, we have $L= 3^{\alpha+2}m(1-G_1)+8\cdot3^k\geq 8\cdot3^k.$

Case 2: Let $0\leq r_1\leq 3$, $r_2= \alpha+2$. Then $G_1=G_2$, $\dfrac{1}{3^{2(\alpha+1)}t^2}\leq G_1\leq \dfrac{1}{3^{2(\alpha+1)}}$ and $G_3=G_4=\dfrac{1}{9}$. Therefore $L= 3^{\alpha+2}m\left(\dfrac{1}{9}-G_1\right)\geq0.$

Case 3: Let $r_1=4$, $0\leq r_2 \leq \alpha+1$. Then $G_2=4G_1$, $\dfrac{1}{4\cdot3^{(\alpha+1)}t^2}\leq G_1\leq \dfrac{1}{4}$, $G_4=4G_3$ and $G_3=\dfrac{1}{4}$. Hence, we have $L=8\cdot3^k.$

Case 4: Let $r_1=4$, $ r_2 =\alpha+2$. Then $G_2=4G_1$, $\dfrac{1}{4\cdot3^{(\alpha+1)}t^2}\leq G_1\leq \dfrac{1}{4\cdot3^{2(\alpha+1)}}$, $G_4=4G_3$ and $G_3=\dfrac{1}{36}$. Therefore $L=0.$

Therefore $D_{\alpha,m,k}(z)$ is holomorphic at every cusp $c/d$. The weight of $D_{\alpha,m,k}(z)$ is $\ell=\dfrac{3^\alpha m-1}{2}+3^k$ which is a positive integer and the associated character is given by
 \begin{align*}
\chi_2(\bullet)=\left(\dfrac{(-1)^\ell   3^{2 \alpha 3^k+3^\alpha \alpha m+3^\alpha m+3^k-1} m^{3^\alpha m+2\cdot 3^k}}{\bullet}\right).
\end{align*}

 Thus, $D_{\alpha,m,k}(z)\in M_\ell\left(\Gamma_0(N),\chi\right)$ where $\ell$, $N$ and $\chi$ are as above. Therefore by Theorem \ref{Ser}, the Fourier coefficients of $D_{\alpha,m,k}(z)$ are almost divisible by $r=3^k$. Due to \eqref{D amk}, this holds for $\overline{a}_{3^{\alpha}m}(n)$ also. This completes the proof of Theorem \ref{a 3k}.

\section{Proof of Theorem \ref{cong1}}\label{pr cong1}

First we recall the following result of Ono and Taguchi \cite{ot} on the nilpotency of Hecke operators. 
\begin{theorem}\cite[Theorem 1.3 (3)]{ot}\label{thmot}
Let $n$ be a nonnegative integer and $k$ be a positive integer. Let $\chi$ be a quadratic Dirichlet character of conductor $9\cdot2^a$. Then there is an integer $c \geq 0$ such that for every $f(z) \in M_k(\Gamma_0(9\cdot2^a),\chi) \cap \mathbb{Z}[[q]]$ and every $t\geq 1$,
\begin{align*}
f(z)|T_{p_1}|T_{p_2}|\cdots|T_{p_{c+t}}\equiv 0\pmod{2^t}
\end{align*}
whenever the primes $p_1, \ldots, p_{c+t}$ are coprime to 6.
\end{theorem}

Now, we apply the above theorem to the modular form $B_{1,1,k}(z)$ to prove Theorem \ref{cong1}.

Putting $\alpha=1$ and $m=1$ in \eqref{B amk}, we find that
\begin{align*}
B_{1,1,k}(z)\equiv\displaystyle\sum_{n=0}^{\infty}\overline{a}_{3}(n)q^{24n}\pmod{2^{k+1}},
\end{align*}
which yields
\begin{align}
B_{1,1,k}(z):=\displaystyle\sum_{n=0}^{\infty}F_k\left(n\right)q^{n}\equiv\displaystyle\sum_{n=0}^{\infty}\overline{a}_{3}\left(\dfrac{n}{24}\right)q^{n}\pmod{2^{k+1}}\label{B 11k}.
\end{align}

Now, $B_{1,1,k}(z)\in M_{2^{k-1}+1}\left(\Gamma_0(9\cdot 2^6),\chi_3\right)$ for $k\geq 1$ where $\chi_3$ is the associated character (which is $\chi_1$ evaluated at $\alpha=1$ and $m=1$). In view of Theorem \ref{thmot}, we find that there is an integer $c\geq 0$ such that for any $d\geq 1$,
\begin{align*}
B_{1,1,k}(z)\mid T_{p_1}\mid T_{p_2}\mid\cdots\mid T_{p_{c+d}}\equiv 0\pmod{2^d}
\end{align*}
whenever $p_1,\dots, p_{c+d}$ are coprime to 6. It follows from the definition of Hecke operators that if $p_1,\dots, p_{c+d}$ are distinct primes and if $n$ is coprime to $p_1\cdots p_{c+d}$,
then
\begin{align}
F_k\left(p_1\cdots p_{c+d}\cdot n\right)\equiv 0\pmod{2^d}.\label{Fa2}
\end{align}

Combining \eqref{B 11k} and \eqref{Fa2}, we complete the proof of the theorem.

\section{Concluding Remarks}\label{co}
\begin{enumerate}
\item[(1)] Theorems \ref{a 2k} and \ref{a 3k} of this paper and Theorem 1.8 of \cite{grs} give us the arithmetic densities of $\overline{a}_t(n)$ for odd $t$ and similar techniques can not be used to obtain the arithmetic density of $\overline{a}_t(n)$ when $t$ is even. It would be interesting to study the arithmetic density of $\overline{a}_t(n)$ for even values of $t$.

\item[(2)] Computational evidence suggests that there are Ramanujan type congruences for $\overline{a}_t(n)$ modulo powers of 2, 3 and other primes $\geq 5$ for various $t$ which are not covered by the results of \cite{bb} and \cite{grs}. We encourage the readers to find new congruences for $\overline{a}_t(n)$. 

\item[(3)] Asymptotic formulas for partition functions and other related functions have been widely studied in the literature. For instance, the asymptotic formulas for $p(n)$ and $c_t(n)$ were obtained by Hardy and Ramanujan \cite{hr} and Anderson \cite{and} respectively. It will be desirable to find an asymptotic formula for $\overline{a}_t(n)$.

\item[(4)] Some relations connecting $\overline{a}_t(n)$ and $c_t(n)$ have been discussed in \cite{bb}. A combinatorial treatment to $\overline{a}_t(n)$ might reveal more interesting partition theoretic connections of $\overline{a}_t(n)$.
\end{enumerate}

\section{Acknowledgement}
The author is extremely grateful to his Ph.D. supervisor, Prof. Nayandeep Deka Baruah, for his guidance and encouragement. The author is indebted to Prof. Rupam Barman for many helpful comments and suggestions. The author was partially supported by the Council of Scientific \& Industrial Research (CSIR), Government of India under the CSIR-JRF scheme (Grant No. 09/0796(12991)/2021-EMR-I). The author thanks the funding agency.

\end{document}